\documentclass[11pt,a4paper]{article}
\usepackage[T1]{fontenc}
\usepackage[utf8]{inputenc}
\usepackage{lmodern,amsmath,amssymb,amsthm,mathtools}
\usepackage[margin=27mm]{geometry}
\usepackage{microtype,enumitem,mathrsfs,needspace}
\usepackage[colorlinks=true,linkcolor=blue,urlcolor=blue,citecolor=blue]{hyperref}
\newtheorem{theorem}{Theorem}[section]
\newtheorem{proposition}[theorem]{Proposition}
\newtheorem{lemma}[theorem]{Lemma}
\newtheorem{corollary}[theorem]{Corollary}
\theoremstyle{definition}
\newtheorem{definition}[theorem]{Definition}
\newtheorem{example}[theorem]{Example}
\theoremstyle{remark}
\newtheorem{remark}[theorem]{Remark}
\newcommand{\C}{\mathbb C}\newcommand{\K}{\mathbb K}\newcommand{\Lop}{\mathbb L}
\newcommand{\Cl}{\mathrm{Cl}_1}\newcommand{\op}{\mathrm{op}}\newcommand{\id}{\mathrm{id}}
\newcommand{\KK}{KK^G}\newcommand{\GK}{\mathsf{GK}^G}
\newcommand{\GKd}{\mathsf{GK}^G_{d,e}}
\newcommand{\gr}{\mathrm{gr}}\newcommand{\eps}{\epsilon}
\newcommand{\Ind}{\operatorname{Index}}\newcommand{\Dom}{\operatorname{Dom}}
\newcommand{\Spec}{\operatorname{spec}}

\title{The universal property of graded \texorpdfstring{\(KK^G\)}{KKG}-theory}
\date{21 October 2026} 	

\author{Bernhard Burgstaller\thanks{%
Departamento de Matematica, 
Universidade Federal de Santa Catarina,  
CEP 88.040-900, 
Florianópolis-SC, 
Brasil. 
\newline
E-mail: \href{mailto:bernhardburgstaller@yahoo.de}{\nolinkurl{bernhardburgstaller@yahoo.de}}}}

\begin{document}
\maketitle

\begingroup
\renewcommand{\thefootnote}{}
\footnotetext{%
}
\endgroup


\begin{abstract}
To each Kasparov cycle $z$ we attach inclusions $B\xrightarrow{j_z}Y_z\xrightarrow{i_z}X_z$ satisfying $[z]j_zi_z=s_zi_z$. A lift of a second factor along $k_z=j_zi_z$ therefore reduces the product to precomposition by a homomorphism. An additional even symmetry gives a retraction of $k_z$; this includes all cycles between trivially graded algebras and the splitting cycles of graded split extensions. The resulting cycle presentation recovers the ungraded equivariant Kasparov category. For arbitrary gradings, explicit tensor and Dirac relations give fusion by passing to a shifted and doubled cycle. We provide concrete Dirac cycles and a compact verification of the two scalar product identities used in this argument. Finally, an index obstruction shows that the corner of the original cycle need not admit a retraction, even after replacing that cycle by an equivalent representative.
\end{abstract}
\tableofcontents
\clearpage

\section*{Introduction}
\addcontentsline{toc}{section}{Introduction}
Kasparov theory for ungraded $C^*$-algebras admits a categorical description by stability, homotopy invariance and split exactness; see Higson \cite{Higson}, Cuntz \cite{Cuntz}, and the equivariant generators-and-relations treatment in \cite{Burgstaller}. The purpose of this note is to organize the fusion of cycles through a different elementary identity, also for graded $C^*$-algebras. To an $(A,B)$-cycle $z=(E,s,F)$ we attach algebra homomorphism inclusions
\[
 B\xrightarrow{j_z}Y_z\xrightarrow{i_z}X_z,\qquad k_z=j_zi_z,
\]
satisfying
\begin{equation}\label{eq1}
 [z]k_z=s i_z
\end{equation}
in Section~\ref{sec:conventions}. Products are written from left to right. Thus multiplication by the auxiliary corner turns the cycle class into the class of its represented homomorphism.

The practical point is a lifting problem, discussed in Section~\ref{sec:sym}. Given $w\in\KK(B,C)$, suppose that a class $w'\in\KK(X_z,C)$ satisfies $k_zw'=w$. Then
\begin{equation}\label{eq:intro-fusion}
 [z]w=[z]k_zw'=s i_z w'.
\end{equation}
Once a cycle for $w'$ is available, the last expression is represented by precomposing its representation with the homomorphism $s i_z$. In this sense, computing the product reduces to providing a lift to the larger algebra $X_z$. If $k_z$ has a retraction $\gamma_z\in\KK(X_z,B)$, a lift always exists: one may take $w'=\gamma_z w$. Constructing a representative of this lift remains the analytic step; for classes of algebras beyond $C^*$-algebras, one may instead seek the required lift directly.

An additional even symmetry, introduced in Definition~\ref{def1}, supplies a retraction explicitly in Theorem~\ref{thm:sym}. This includes every cycle between trivially graded algebras, as well as the usual splitting cycles for graded split extensions 
recalled in Example \ref{prop:split}. Section~\ref{sec5} gives the corresponding cycle presentation and proves fusion at such cycles in Proposition~\ref{prop:fuse}, yielding a faithful subcategory in Corollary~\ref{cor:subcategory}. Corollary~\ref{cor:ungraded} recovers the usual ungraded equivariant Kasparov category. Its practical benefit is that the equivariant actions remain on the Hilbert modules throughout: the fusion calculation requires no explicit transport between the actions belonging to the two represented summands, as in split-extension presentations such as \cite{Burgstaller}.

For arbitrary graded cycles, Section~\ref{sec:dirac} adds specified tensor compatibilities and two Dirac relations. A shifted and doubled cycle has the required symmetry; the Dirac relations, whose scalar identities are proved in Lemma~\ref{lem:ed}, recover the original class and remove the extra branch. Theorem~\ref{thm:dirac-fusion} proves that the enlarged presentation is the full graded equivariant Kasparov category. Sections~\ref{sec:obstruction} and~\ref{sec:invariant} explain why this enlargement is useful: an explicit index obstruction excludes a retraction at the corner of the original cycle for every representative of a certain class. The new fusion argument uses a different corner and is consistent with this obstruction. This corrects the earlier version of the present paper \cite{flawedpaper}, whose graded argument did not justify the required retraction.

Alternative descriptions of graded Kasparov theory are provided by Haag's algebraic and homomorphism pictures \cite{Haag,HaagExt} and their equivariant development by Meyer \cite{MeyerEquivariant}. In particular, Meyer's treatment relates graded theory to ungraded equivariant theory with an additional $\mathbb Z_2$-action and describes the product in that setting. The presentation here is organized around \eqref{eq1}, the lifting construction of Sections~\ref{sec:sym} and~\ref{sec5}, and the explicit Dirac relations of Section~\ref{sec:dirac}. For the broader categorical viewpoint, see \cite{Meyer}.

\Needspace{9\baselineskip}
\subsection*{AI Disclosure}
\addcontentsline{toc}{subsection}{AI Disclosure}
This is joint work with OpenAI ChatGPT Astra 5.1 or 6. The positive results in Sections~\ref{sec:conventions}--\ref{sec5} are by the human author and form part of the correction of the flawed paper \cite{flawedpaper}. The negative results (obstructions) in Sections~\ref{sec:obstruction} and~\ref{sec:invariant} were found by ChatGPT. The general formula \eqref{eq:form} was found by Google Gemini Pro 3.1; see Remark~\ref{remark}. The argument in Section~\ref{sec:conventions} that $X_z$ is norm-closed for trivially graded $A$, including the formulas for $F_i$, is by the human author; its extension to arbitrary graded $A$ is by ChatGPT.

Section~\ref{sec:dirac} was mainly done by ChatGPT. The human author observed from the beginning that one could axiomatize $z+z^\sharp$, but not the desired doubling $2z$. After the counterexamples in Sections~\ref{sec:obstruction} and~\ref{sec:invariant} had been found, ChatGPT wrote a larger working note for the human author in another context and briefly introduced elements $d$ and $e$ that eliminate the superfluous summand $z^\sharp$. The human author immediately recognized that these elements 
supply the missing fusion argument by adding the Dirac relations (D), and instructed ChatGPT to work this out. The resulting construction in Section~\ref{sec:dirac} uses the retraction $\gamma_u$ of a shifted and doubled cycle and then recovers the original class through the Dirac relations.

The text was written by ChatGPT; the AI disclosure is by the human author. 


\section{Conventions and the algebras attached to a cycle}\label{sec:conventions}
All algebras and Hilbert modules are $\mathbb Z_2$-graded. Algebras are separable and modules countably generated. We work in an equivariant Kasparov setting with products and graded Morita equivalence. This includes second-countable locally compact Hausdorff groupoids $G$, with the $C_0(G^{(0)})$-linearity and continuity requirements of Le Gall \cite{LeGall}; for groups see \cite{Kasparov88}. The algebras below carry their induced continuous actions; the full adjointable-operator algebra is not asserted to be a continuous $G$-algebra. The counterexample is non-equivariant.

An $(A,B)$-cycle is written $z=(E,\pi,F)$, with the module first. Representations are even, operators $F$ are odd, and commutators are graded. Products of morphisms are written from left to right: $xy$ means first $x$, then $y$. The grading automorphism of an algebra $A$ is $\eps_A$. The cycle conditions are
\[
 [\pi(a),F]_{\gr},\quad \pi(a)(F-F^*),\quad
 \pi(a)(F^2-1)\in\K_B(E)\qquad(a\in A),
\]
together with local equivariance of $F$ modulo compact operators. In groupoid notation the latter means
\[
 r^*\pi(b)\bigl(U(s^*F)U^{-1}-r^*F\bigr)
 \in r^*\K_B(E)\qquad(b\in r^*A),
\]
where $s,r$ are source and range and $U:s^*E\to r^*E$ is the module action. We use the standard Kasparov homotopy relation and product \cite{Kasparov81,Kasparov88,LeGall}. For $E^{\op}$ we reverse the module grading and retain its right action and inner product; grading twists on left representations are written explicitly.

Adjoin the zero-representation cycle
\[
 (B,0,0),\qquad \mathcal H_z=E\oplus B,\quad
 \widetilde\pi=\pi\oplus0,\quad \widetilde F=F\oplus0.
\]
We suppress tildes below. Set
\begin{equation}\label{eq:XY}
 J_z=\K_B(\mathcal H_z),\qquad
 Y_z=J_z+\pi(A),\qquad
 X_z=J_z+\pi(A)+F\pi(A).
\end{equation}
These algebraic sums are already closed. Here is a quotient argument that also treats nontrivial gradings. Let
\[
 q:\Lop_B(\mathcal H_z)\longrightarrow\Lop_B(\mathcal H_z)/J_z,
 \qquad \rho=q\pi,\quad f=q(F),\quad\alpha=\eps_A.
\]
The Kasparov relations give
\[
 f\rho(a)=\rho(\alpha(a))f,\qquad
 (f^2-1)\rho(a)=0,\qquad(f-f^*)\rho(a)=0.
\]
The closed $C^*$-subalgebra
\[
 C_\alpha=
 \left\{\begin{pmatrix}a&b\\\alpha(b)&\alpha(a)\end{pmatrix}:a,b\in A\right\}
 \subseteq M_2(A)
\]
therefore has a $*$-homomorphism
\[
 \Phi:C_\alpha\longrightarrow\Lop_B(\mathcal H_z)/J_z,
 \qquad
 \begin{pmatrix}a&b\\\alpha(b)&\alpha(a)\end{pmatrix}
 \longmapsto\rho(a)+\rho(b)f.
\]
Indeed, the multiplication rule on the pairs is
$(a,b)(c,d)=(ac+b\alpha(d),ad+b\alpha(c))$, and the displayed quotient relations give exactly this rule and its adjoint rule. The images of $\rho$ and $\Phi$ are closed, so
\[
 Y_z=q^{-1}(\rho(A)),\qquad X_z=q^{-1}(\Phi(C_\alpha))
\]
are closed. When $A$ is trivially graded, put $F_i=(1+(-1)^iF)/2$ for $i=1,2$. Then
\[
 \pi(A)+F\pi(A)=F_1\pi(A)+F_2\pi(A).
\]
The maps $a\mapsto q(F_i\pi(a))$ are $*$-homomorphisms with orthogonal images, so their sum is closed. Taking the inverse image under $q$ again proves that $X_z$ is closed.

The algebras in \eqref{eq:XY} are graded, and $F$ is a multiplier of $X_z$. For example, for homogeneous $a$,
\[
 \pi(a)F\equiv(-1)^{|a|}F\pi(a)\pmod{J_z},\qquad
 F^2\pi(a)\equiv\pi(a)\pmod{J_z},
\]
which give $FX_z,X_zF\subseteq X_z$. Separability follows from that of $J_z$ and $A$. Equivariance and continuity follow from the local equivariance defects of the cycle; for groupoids the same statements use the corresponding pullbacks.

There are even equivariant homomorphisms
\[
 j_z:B\longrightarrow Y_z,\qquad i_z:Y_z\hookrightarrow X_z,
 \qquad k_z=j_zi_z,
\]
where $j_z$ is left multiplication in the added $B$-corner. The full module $\mathcal H_z$ implements the Morita equivalence between $J_z$ and $B$. Also
\[
 \delta_z=[\mathcal H_z,\id_{Y_z},F]\in\KK(Y_z,B),
 \qquad [z]=\pi\delta_z.
\]
Here $\id_{Y_z}$ denotes its represented inclusion into $\Lop_B(\mathcal H_z)$.

\begin{lemma}[The analytic KK identity]\label{lem:analytic}
For every cycle,
\[
 j_z\delta_z=1_B,\qquad \delta_zj_zi_z=i_z.
\]
\end{lemma}
\begin{proof}
The first identity is restriction to the added $B$-corner, on which $F=0$. For the second, first transport $\delta_z$ along the Morita corner $B\to J_z$. This gives the cycle $(J_z,\id_{Y_z},F)$ with $J_z$ regarded as a right $J_z$-module. Pushing forward to $X_z$ gives the same ideal as a right $X_z$-module.

More generally, if $J\triangleleft X\subseteq M(J)$, $s:A\to X$, and $T\in M(X)$ defines $(J,s,T)$, the right $X[0,1]$-module
\[
 \mathcal V=\{f\in C([0,1],X):f(0)\in J\}
\]
with pointwise $s$ and $T$ gives a homotopy from $(J,s,T)$ to $(X,s,T)$. Its compact operators identify with the ideal $\mathcal V$; the constant defect functions have values in $J$, so belong to this ideal. At the second endpoint, $s(a)\in X=\K_X(X)$, hence scaling $T$ to zero is a Kasparov homotopy. The result is the homomorphism class $[s]$.

In the present case take $J=J_z$, $X=X_z$, $s=i_z$ and $T=F$. The transported cycle above represents $\delta_zj_zi_z$; the homotopy enlarges its right module from $J_z$ to $X_z$ and then scales $F$ to zero, leaving $(X_z,i_z,0)$, which represents $i_z$. The pointwise construction respects the equivariant structure.
\end{proof}

\section{An additional even symmetry}\label{sec:sym}
\begin{definition}\label{def1}
A cycle $z$ as in Section~\ref{sec:conventions} is called \emph{symmetry-admitting} if there is an even equivariant selfadjoint unitary $p\in\Lop_B(E)$ such that
\begin{equation}\label{eq:p}
 [p,\pi(a)]=0\quad(a\in A),\qquad pF+Fp=0.
\end{equation}
\end{definition}
Such a $p$ is equivalent to a decomposition into graded equivariant Hilbert $B$-submodules
\[
 E=E_1\oplus E_2,\quad
 \pi=\begin{pmatrix}\pi_1&0\\0&\pi_2\end{pmatrix},\quad
 F=\begin{pmatrix}0&F_2\\F_1&0\end{pmatrix}.
\]
Indeed $p=\operatorname{diag}(1,-1)$, and conversely $E_{1,2}$ are the ranges of the even projections $(1\pm p)/2$. Each $E_i$ has its own grading. This decomposition is additional to the intrinsic module grading.

Extend $p$ to $E\oplus B$ by $p\oplus1_B$ and set
\begin{equation}\label{eq:form}
 T=-ipF.
\end{equation}
\begin{theorem}[Symmetry lemma]\label{thm:sym}
For a symmetry-admitting cycle, the class
\[
 \gamma_z=[\mathcal H_z,\id_{X_z},-ipF]\in\KK(X_z,B)
\]
is well defined and satisfies
\begin{equation}\label{eq:gamma}
 k_z\gamma_z=1_B.
\end{equation}
\end{theorem}
\begin{proof}
This is a straightforward computation using the element \(\delta_z\) from Section \ref{sec:conventions}. On the distinguished $B$-corner one has $F=T=0$; restriction along $k_z$ therefore gives $(B,\id_B,0)$ together with zero-representation summands, proving \eqref{eq:gamma}.
\end{proof}

\begin{proposition}[Symmetry modulo compacts]\label{prop:mod}
The conclusion remains valid if $p$ is still even, equivariant and a selfadjoint unitary, but \eqref{eq:p} is replaced by
\[
 [p,\pi(a)]\in J_z\quad(a\in A),\qquad pF+Fp\in J_z.
\]
\end{proposition}

\begin{example}[Graded splitting cycles]\label{prop:split}
Let $0\to I\xrightarrow{\iota}D\xrightarrow{r}A\to0$ have an even equivariant section $s:A\to D$, so $r\circ s=\id_A$. Write $\lambda:D\to\Lop_I(I)$ for left multiplication. The standard splitting cycle \cite[Example~17.1.2(b)]{Blackadar}, written in our opposite-grading convention, is
\[
 d_s=\left[I\oplus I^{\op},\ 
 d\longmapsto\begin{pmatrix}\lambda(d)&0\\0&\lambda(s(r(\eps_D(d))))\end{pmatrix},\ 
 \begin{pmatrix}0&1\\1&0\end{pmatrix}\right].
\]
It has the required diagonal/off-diagonal form, with $p=\operatorname{diag}(1,-1)$. Thus Theorem~\ref{thm:sym} applies, with no restriction on the gradings of $I,D,A$.
\end{example}

\begin{remark}[Ungraded cycles and the formula for the retraction]\label{remark}
If $A$ and $B$ are trivially graded, the intrinsic module grading is a $B$-linear even equivariant selfadjoint unitary. It commutes with $\pi(A)$ and anticommutes with $F$, so every such cycle admits the symmetry. The construction retains the original module action and does not require identifying the two graded summands and transporting the action between them.

In a normalized ungraded picture one may arrange $F_1=F_2=1$ and, after adjoining a zero-class cycle, make the representation faithful. The author originally obtained the retraction in this picture by a composition of inverse corner embeddings and splitting morphisms, one associated with unitization and one with the representation. This links the construction to the usual stability and split exactness presentation. The author also knew \eqref{eq:form} in the special case 
$F_1=F_2$; Google Gemini Pro 3.1 independently suggested the general formula $T=-ipF$.
\end{remark}

\begin{lemma}[Universal lifting criterion]\label{thm:lift}
For a fixed cycle $z$, the following are equivalent in $\KK$:
\begin{enumerate}[label=(\roman*)]
\item There exists $\gamma_z\in\KK(X_z,B)$ with $k_z\gamma_z=1_B$.
\item For every target $D$ and every $w\in\KK(B,D)$ there is $y\in\KK(X_z,D)$ with $k_zy=w$.
\end{enumerate}
\end{lemma}
\begin{proof}
Use $y=\gamma_zw$ for (i)$\Rightarrow$(ii), and $D=B$, $w=1_B$ for the converse.
\end{proof}
Whenever this holds, Lemma~\ref{lem:analytic} also implies
\begin{equation}\label{eq2}
 \delta_z=\delta_zk_z\gamma_z=i_z\gamma_z.
\end{equation}
\section{The cycle presentation and fusion}\label{sec5}
For each padded Kasparov $(A,B)$-cycle $z=(\mathcal H_z,s,F_z)$, retain the algebras and maps $Y_z,X_z,j_z,i_z$ of Section~\ref{sec:conventions}. Introduce one synthetic arrow
\[
 \Delta_z:Y_z\longrightarrow B.
\]
Together with graded equivariant homomorphisms, these generate a free additive category $\mathsf{GK}^G_{\mathrm{free}}$ on the graded $G$-algebras. A suitably small collection of cycles is understood.

We call an expression
\[
 a\Delta_z:C\longrightarrow B,
 \qquad a:C\longrightarrow Y_z\text{ a homomorphism},
\]
a \emph{simple form}. The map $a$ need not be the original representation $s$ of $z$. In particular, the generator itself is the simple form $\id_{Y_z}\Delta_z$.

Define $\GK$ to be the quotient by the composition and identity relations for homomorphisms and the following relations. For categorical background see \cite{Meyer}.
\begin{enumerate}[label=\textnormal{\arabic*.},leftmargin=*]
\item \emph{Cycle invariance.} Set $a\Delta_z=b\Delta_v$ whenever the associated Kasparov $(C,B)$-cycles are equivalent under the standard Kasparov homotopy relation and even equivariant unitary equivalence.
\item \emph{Additivity.} A simple form associated with a direct sum of two cycles equals the sum of their simple forms. A simple form for the inverse class equals the negative.
\item \emph{Homomorphism normalization.} For a homomorphism $h:C\to B$, let $z_h$ be its padded homomorphism cycle
\[
 z_h=(B\oplus B,s_h,0),\qquad s_h(c)=\operatorname{diag}(h(c),0).
\]
Then $s_h\Delta_{z_h}=h$. The second copy of $B$ is the distinguished zero-representation corner.
\item \emph{The KK relation.} For every $z$,
\[
 \Delta_zj_zi_z=i_z.
\]
\end{enumerate}
Homotopy invariance for homomorphisms follows from cycle invariance and homomorphism normalization.

\begin{lemma}[Comparison on simple forms]\label{lem:simple}
There is an additive functor $\mathcal A:\GK\to\KK$ defined by
\[
 \mathcal A(h)=[h],\qquad \mathcal A(\Delta_z)=\delta_z.
\]
Choosing any simple representative defines an additive assignment
\[
 \mathcal B:\KK(A,B)\to\GK(A,B),\qquad
 \mathcal B([\mathcal H_z,s,F_z])=s\Delta_z.
\]
One has
\[
 \mathcal A\mathcal B=\id,\qquad
 \mathcal B\mathcal A(a\Delta_z)=a\Delta_z.
\]
Moreover, for every homomorphism $f:A'\to A$,
\begin{equation}\label{eq:precomp}
 f\mathcal B(x)=\mathcal B(fx).
\end{equation}
\end{lemma}
\begin{proof}
All relations hold in $\KK$, the last by Lemma~\ref{lem:analytic}. Cycle invariance makes $\mathcal B$ well defined, and additivity makes it additive. The inverse identities on simple forms follow from
\[
 \mathcal A(a\Delta_z)=[\mathcal H_z,a,F_z].
\]
Finally $f(s\Delta_z)=(s\circ f)\Delta_z$ by associativity and homomorphism composition. Its associated cycle is the pullback along $f$, proving \eqref{eq:precomp}. In particular $\mathcal B(\delta_z)=\Delta_z$. No multiplicativity of $\mathcal B$ is assumed.
\end{proof}

\begin{proposition}[Fusion at a retractable cycle]\label{prop:fuse}
If $k_z\gamma_z=1_B$ in $\KK$, then for every $w\in\KK(B,D)$,
\[
 \Delta_z\mathcal B(w)=\mathcal B(i_z\gamma_z w).
\]
The right-hand side is a simple form.
\end{proposition}
\begin{proof}
Put $y=\gamma_z w$. Equation~\eqref{eq:precomp} gives
$k_z\mathcal B(y)=\mathcal B(k_z y)=\mathcal B(w)$. Multiplication by $\Delta_z$ and the KK relation give
\[
 \Delta_z\mathcal B(w)=\Delta_zk_z\mathcal B(y)
 =i_z\mathcal B(y)=\mathcal B(i_zy).
\]
Only precomposition by homomorphisms has been used.
\end{proof}

\begin{corollary}\label{cor:subcategory}
Let $\mathscr S$ be the additive subcategory of $\GK$ generated by homomorphisms and simple forms $\pi\Delta_z$ for symmetry-admitting cycles, also allowing Proposition~\ref{prop:mod}. Then $\mathcal A|_{\mathscr S}$ is faithful and is an isomorphism onto the additive subcategory of $\KK$ generated by their images, with inverse the restriction of $\mathcal B$. Every morphism in $\mathscr S$ has a simple form.
\end{corollary}
\begin{proof}
Reduce a composable word from right to left. Its rightmost factor is a simple form, using homomorphism normalization if needed. Each additional left homomorphism is handled by \eqref{eq:precomp}; each additional left factor $\pi\Delta_z$ is handled by Proposition~\ref{prop:fuse} and then \eqref{eq:precomp}. At each stage the left factor is one of the original symmetry-admitting factors; the intermediate product need not admit a symmetry. Additivity handles sums and negatives. Hence every arrow in $\mathscr S$ equals $\mathcal B$ of its image, proving faithfulness and the image assertion.
\end{proof}

\begin{corollary}[The ungraded cycle presentation]\label{cor:ungraded}
Restrict $\mathscr S$ to its trivially graded objects. Then
\[
 \mathcal A:\mathscr S_{\mathrm{triv}}\xrightarrow{\ \cong\ }
 KK^G_{\mathrm{triv}}
\]
is an isomorphism of additive categories, with inverse $\mathcal B$.
\end{corollary}
\begin{proof}
For trivially graded endpoints the intrinsic module grading supplies the symmetry of Remark~\ref{remark}. Thus every $KK^G$-class lies in the image of Corollary~\ref{cor:subcategory}, and its faithfulness proves the isomorphism.
\end{proof}
This recovers the usual ungraded equivariant Kasparov category, hence its homotopy invariance, stability and split exactness. The graded splitting morphisms of Example~\ref{prop:split} also belong to the faithful subcategory. Remark~\ref{remark} and \eqref{eq2} explain the connection with the usual corner and splitting generators. The classical universal characterization is the background theorem \cite{Higson,Burgstaller}; the present argument establishes the stated cycle presentation. Auxiliary algebras $Y_z,X_z$ may be graded even when the endpoints are ungraded.

\section{Dirac relations and fusion for arbitrary graded cycles}\label{sec:dirac}
We retain the cycle presentation of Section~\ref{sec5}. In particular,
\begin{equation}\label{eq:dirac-basic}
 \mathcal A\mathcal B=\id,\qquad \mathcal B([z])=s_z\Delta_z,
 \qquad h\mathcal B(x)=\mathcal B(hx)
\end{equation}
for homomorphisms $h$. Multiplicativity of $\mathcal B$ is not assumed. We will use a corner attached to a shifted and doubled cycle and remove the extra branch by two fixed Dirac relations.

\subsection{The fixed cycles and their two products}\label{subsec:ed}
Put $S=C_0(\mathbb R)$ with reflection grading $\alpha(f)(t)=f(-t)$, and $TA=S\widehat\otimes A$. For a cycle $z=(E,\pi,F)$ write
\[
 z^\sharp=(E,\pi,-F),
\]
with unchanged module grading and representation. This operation descends to classes and does not denote the additive inverse. Reversing the signs of both input operators and a product operator preserves the connection and positivity conditions defining a Kasparov product. Consequently the analytic exterior product satisfies
\begin{equation}\label{eq:sharp-exterior}
 (a\boxtimes x)^\sharp=a^\sharp\boxtimes x^\sharp.
\end{equation}
Only exterior products with a non-equivariant scalar factor are used below. This factor carries the trivial action; for groupoids the $C_0(G^{(0)})$-structure and action remain on the other factor.

We now specify all the Dirac data. On the Hilbert $S$-module $S$, with its reflection grading, let
\[
 (Q_S f)(t)=tf(t),\qquad
 \Dom Q_S=\{f\in C_0(\mathbb R):tf(t)\in C_0(\mathbb R)\}.
\]
On $H=L^2(\mathbb R)$ use the grading $R\xi(t)=\xi(-t)$, the multiplication representation $m(f)=M_f$, and
\begin{equation}\label{eq:D-domain}
 D=-R\partial=-iRP,\qquad P=-i\partial,\qquad
 \Dom D=H^1(\mathbb R).
\end{equation}
Here $\partial\xi=\xi'$ is the \emph{weak, or distributional, derivative}, and
$H^1(\mathbb R)=\{\xi\in L^2(\mathbb R):\xi'\in L^2(\mathbb R)\}$.
Thus $(D\xi)(t)=-\xi'(-t)$. For a selfadjoint regular operator $L$, write $F_L=L(1+L^2)^{-1/2}$ for its bounded transform. Our fixed classes are
\begin{equation}\label{eq:ed-cycles}
 e=[S,\mathrm{scal},F_{Q_S}],\qquad
 d=[H,m,F_D],\qquad
 d^\sharp=[H,m,-F_D].
\end{equation}
Their types are $e\in KK(\C,S)$ and $d,d^\sharp\in KK(S,\C)$.

\Needspace{11\baselineskip}
\begin{lemma}[The two scalar Dirac products]\label{lem:ed}
For the cycles \eqref{eq:ed-cycles},
\begin{equation}\label{eq:ed}
 ed=1_\C,\qquad ed^\sharp=0.
\end{equation}
\end{lemma}
\begin{proof}
We use Kucerovsky's unbounded product criterion \cite[Theorem~13]{Kucerovsky}, in the formulation of \cite[Theorem~2.5]{vandenDungen}. The proposed product operators will be $L_+=Q+D$ and $L_-=Q-D$ on $H$. Besides showing that these are unbounded Kasparov cycles, we must check three conditions: a bounded connection commutator comparing $L_\varepsilon$ with the second-factor operator $\varepsilon D$; inclusion of the product domain in the domain of the first-factor operator $Q$; and a lower bound for their mixed quadratic form. These checks identify the products; the kernels then determine their classes.

First, the input data in \eqref{eq:ed-cycles} are valid cycles. The coordinate multiplier $Q_S$ is odd, regular and selfadjoint, with $(1+Q_S^2)^{-1}\in S$. The operator $D$ is odd and selfadjoint on $H^1$: use $\partial^*=-\partial$ and $R\partial=-\partial R$. Moreover $D^2=-\partial^2$. For homogeneous $g\in C_c^\infty(\mathbb R)$ of parity $p$,
\begin{equation}\label{eq:D-comm}
 [D,M_g]_{\gr}=(-1)^pM_{g'}R.
\end{equation}
Local compactness follows from the compact inclusion of $H^1$ into $L^2$ on bounded intervals, first for compactly supported $g$ and then by approximation for $g\in C_0(\mathbb R)$. Thus the displayed data give unbounded Kasparov cycles.

The even unitary $S\widehat\otimes_S H\to H$, $f\otimes\xi\mapsto f\xi$, identifies $Q_S\widehat\otimes1$ with coordinate multiplication $Q$ on $H$. For $\varepsilon\in\{+1,-1\}$ put
\[
 L_\varepsilon=Q+\varepsilon D,\qquad
 V=\Dom L_\varepsilon=H^1(\mathbb R)\cap\Dom Q,
 \qquad\Dom Q=\{\xi:t\xi(t)\in L^2(\mathbb R)\}.
\]
We first justify this choice of domain. On the Schwartz space $\mathscr S(\mathbb R)$, the ordinary differentiation rule and $RQ=-QR$ give $QD+DQ=-R$. Taking its quadratic form yields
\[
 2\operatorname{Re}\langle Q\xi,D\xi\rangle
 =-\langle\xi,R\xi\rangle.
\]
This identity controls the cross term in the square of $Q+\varepsilon D$. We use the inner product linear in its second argument. The space $C_c^\infty(\mathbb R)$, and therefore also $\mathscr S(\mathbb R)$, is dense in $V$ for the complete norm
\[
 \|\xi\|_V=\bigl(\|\xi\|^2+\|\xi'\|^2+\|Q\xi\|^2\bigr)^{1/2}.
\]
Indeed, multiplication by smooth cutoffs first approximates all three terms of this norm; smoothing the compactly supported functions then gives the required approximation. Both sides of the form identity are continuous in this norm, so the identity extends to every $\xi\in V$. This uses only $Q\xi$ and $D\xi$, and does not require the separate expressions $QD\xi$ or $DQ\xi$ to exist. Expanding the square now gives
\begin{equation}\label{eq:dirac-graph}
 \|L_\varepsilon\xi\|^2
 =\|Q\xi\|^2+\|\xi'\|^2-\varepsilon\langle\xi,R\xi\rangle.
\end{equation}
Since $\|R\|=1$, the graph norm of $L_\varepsilon$ is equivalent to $\|\cdot\|_V$. Thus the proposed operator is closed on $V$; symmetry follows from that of $Q$ and $D$. This establishes a closed operator on a precise domain before the spectral calculation below.

For completeness, selfadjointness and the spectrum are visible in the Hermite basis. Let $h_n$ be the normalized Hermite functions, with
$Rh_n=(-1)^nh_n$ and $(-\partial^2+Q^2)h_n=(2n+1)h_n$. Their usual differentiation and multiplication recursions give
\begin{align*}
 L_+h_0&=0,&
 L_+|_{\operatorname{span}\{h_{2k-1},h_{2k}\}}
 &=2\sqrt{k}\begin{pmatrix}0&1\\1&0\end{pmatrix}\quad(k\geq1),\\
 &&L_-|_{\operatorname{span}\{h_{2k},h_{2k+1}\}}
 &=\sqrt{4k+2}\begin{pmatrix}0&1\\1&0\end{pmatrix}\quad(k\geq0).
\end{align*}
The maximal orthogonal sums of these finite selfadjoint blocks are selfadjoint and have the finite Hermite span as a core. Each is contained in the closed symmetric operator $L_\varepsilon$ just constructed; selfadjointness then forces equality. Thus its domain is exactly $V$. The block eigenvalues tend to infinity in modulus, so both operators have compact resolvent and define unbounded $(\C,\C)$-cycles.

We can now verify the three product conditions announced above. For the connection condition, use the graded subspace $C_c^\infty(\mathbb R)\subset\Dom Q_S$ of the first module $S$. It is dense even in the graph norm $\|g\|_\infty+\|tg\|_\infty$, by cutoff and smoothing. Every test vector is a sum of an even and an odd one, so it suffices to take homogeneous $g$ of parity $p$. Under the tensor-product identification, creation by $g$, namely $\xi\mapsto g\otimes\xi$, becomes $T_g=M_g$, with adjoint $T_g^*=M_{\bar g}$. The weak product rule $(g\xi)'=g'\xi+g\xi'$ gives
\[
 M_g:H^1\longrightarrow V,\qquad M_{\bar g}:V\longrightarrow H^1.
\]
These maps preserve the domains needed to form the graded commutator of $\operatorname{diag}(L_\varepsilon,\varepsilon D)$ with the off-diagonal matrix having entries $T_g$ and $T_g^*$. Its upper entry is the bounded operator
\begin{equation}\label{eq:dirac-connection}
 A_{\varepsilon,g}
 =L_\varepsilon M_g-(-1)^p M_g\varepsilon D
 =M_{tg}+\varepsilon(-1)^pM_{g'}R;
\end{equation}
because both $tg$ and $g'$ are bounded. The lower entry is $-(-1)^pA_{\varepsilon,g}^*$ and is bounded as well. This proves the connection condition, including preservation of domains.

The domain condition asks that $\Dom L_\varepsilon$ be contained in $\Dom(Q_S\widehat\otimes1)$; under our identification this is exactly the inclusion $V\subseteq\Dom Q$ already built into the definition. Finally, positivity means that the mixed form of $Q$ and $L_\varepsilon$ has a uniform lower bound on that domain. The form identity used in \eqref{eq:dirac-graph} gives, for every $\xi\in V$,
\begin{equation}\label{eq:dirac-positive}
 2\operatorname{Re}\langle Q\xi,L_\varepsilon\xi\rangle
 =2\|Q\xi\|^2-\varepsilon\langle\xi,R\xi\rangle
 \geq-\|\xi\|^2.
\end{equation}
Thus the positivity constant can be chosen as $-1$ for both signs. The earlier extension of the form identity is what makes this estimate valid on the whole product domain. The criterion therefore identifies $[H,\mathrm{scal},F_{L_+}]$ with $ed$ and $[H,\mathrm{scal},F_{L_-}]$ with $ed^\sharp$.

It remains to read off these classes. The blocks give
\[
 \Spec(L_+^2)=\{0,4,8,\ldots\},\qquad
 \Spec(L_-^2)=\{2,6,10,\ldots\}.
\]
The first kernel is the single even Gaussian $h_0=\pi^{-1/4}e^{-t^2/2}$; the second kernel is zero. The function $\operatorname{sgn}_0$, defined to be zero at zero, is continuous on $\Spec(L_\varepsilon)$ because of the spectral gap. Equivalently, a continuous odd function on $\mathbb R$ that equals $\pm1$ for $\pm t\geq1$ and equals $t$ on $[-1,1]$ has precisely these spectral values. Continuous functional calculus therefore gives $G_\varepsilon=\operatorname{sgn}_0(L_\varepsilon)$; no measurable functional calculus is needed.

On a nonzero block with modulus $\lambda$, the difference $G_\varepsilon-F_{L_\varepsilon}$ has norm $1-\lambda/\sqrt{1+\lambda^2}\to0$. It is consequently compact, so the straight line between these two odd selfadjoint operators is a Kasparov homotopy. Off the kernel, $G_\varepsilon$ is an odd selfadjoint unitary and the scalar cycle is degenerate. The remaining even line for $L_+$ represents $1_\C$, while $L_-$ leaves no summand. This proves both identities in \eqref{eq:ed}.
\end{proof}

\subsection{The additional relations}
The construction adds relations to the original presentation, with unchanged objects and arrow generators. Set
\begin{align*}
 E_A&=\mathcal B(e\boxtimes1_A):A\longrightarrow TA,\\
 D_A&=\mathcal B(d\boxtimes1_A):TA\longrightarrow A,\qquad
 D_A^-=\mathcal B(d^\sharp\boxtimes1_A):TA\longrightarrow A.
\end{align*}
These are names for existing simple arrows. Use the following notation and relation schemes:
\begin{enumerate}[label=\textnormal{(T\arabic*)},leftmargin=*]
\item \emph{Tensor notation.} For homomorphisms $h$ and classes $x\in\KK(A,B)$, write
\[
 T(h)=\id_S\widehat\otimes h,\qquad
 T(\mathcal B(x))=\mathcal B(1_S\boxtimes x).
\]
For $x=[h]$ the two conventions agree by homomorphism normalization. This is notation only; no operation on arbitrary words or functoriality of $T$ is assumed.
\item \emph{Scalar-cycle relation.} For every $x\in\KK(A,B)$ impose
\begin{equation}\label{eq:T2}
 E_A\mathcal B(1_S\boxtimes x)=\mathcal B(x)E_B.
\end{equation}
\item \emph{Compatibility with the fixed Dirac exterior products.} For every such $x$ impose
\begin{equation}\label{eq:T3}
 \begin{aligned}
 \mathcal B(1_S\boxtimes x)D_B&=\mathcal B(d\boxtimes x),\\
 \mathcal B(1_S\boxtimes x)D_B^-&=\mathcal B(d^\sharp\boxtimes x).
 \end{aligned}
\end{equation}
The right-hand sides are simple arrows of chosen cycles for these analytic exterior-product classes; cycle invariance makes the choices immaterial.
\end{enumerate}
In addition impose the \emph{Dirac relations}, with all tensor copies:
\begin{equation*}\tag{D}\label{eq:D-relations}
 E_AD_A=1_A,\qquad E_AD_A^-=0\qquad\text{for every object }A.
\end{equation*}
Thus (T1) only specifies notation, whereas (T2), (T3) and (D) are additional relations. Every term uses arrows already present in $\GK$.

Let $\GKd$ be the quotient by these relations. The comparison factors through an additive functor
\[
 \mathcal A:\GKd\longrightarrow\KK.
\]
Indeed, (T2) and (T3) hold by analytic exterior-product compatibility, and (D) follows from Lemma~\ref{lem:ed}. The assignment $\mathcal B$ is followed by passage to the quotient; it remains additive and satisfies \eqref{eq:dirac-basic}.

\begin{lemma}[Contraction of the two Dirac branches]\label{lem:contraction}
For every $x\in\KK(A,B)$, in $\GKd$ one has
\begin{equation}\label{eq:contraction}
 E_A\mathcal B(d\boxtimes x)=\mathcal B(x),\qquad
 E_A\mathcal B(d^\sharp\boxtimes x)=0.
\end{equation}
\end{lemma}
\begin{proof}
Using (T3), (T2) and (D), respectively, gives
\[
 E_A\mathcal B(d\boxtimes x)
 =E_AT(\mathcal B(x))D_B
 =\mathcal B(x)E_BD_B=\mathcal B(x).
\]
Replacing $D_B$ by $D_B^-$ proves the second identity. It holds for every $x$, without any condition involving operator sign reversal.
\end{proof}

\subsection{The retraction for the shifted and doubled cycle}
Let $z\in\KK(A,B)$ and choose a selfadjoint cycle $w=(L,\pi,K)$ representing $d\boxtimes z$. Form
\begin{equation}\label{eq:doubled}
 u=w\oplus w^\sharp\oplus(B,0,0),\qquad
 \mathcal H_u=L\oplus L\oplus B.
\end{equation}
The two copies of $L$ have the same grading and action. The last summand is the single distinguished padding summand. For $a\in TA$, the representation and operator are
\[
 \Pi(a)=\begin{pmatrix}\pi(a)&0&0\\0&\pi(a)&0\\0&0&0\end{pmatrix},
 \qquad F=\begin{pmatrix}K&0&0\\0&-K&0\\0&0&0\end{pmatrix}.
\]
Construct $J_u,Y_u,X_u,s_u,j_u,i_u,k_u$ as in Section~\ref{sec:conventions}. In particular $s_u:TA\to Y_u$ is $\Pi$, and $k_u:B\to X_u$ is the last diagonal corner. Set
\begin{equation}\label{eq:doubled-symmetry}
 p=\begin{pmatrix}0&1&0\\1&0&0\\0&0&1\end{pmatrix},\qquad
 T_u=-ipF=\begin{pmatrix}0&iK&0\\-iK&0&0\\0&0&0\end{pmatrix}.
\end{equation}
The operator $p$ is an even equivariant selfadjoint unitary, commutes with $\Pi(TA)$ and anticommutes with $F$. Theorem~\ref{thm:sym} and \eqref{eq2} therefore give
\begin{equation}\label{eq:doubled-retraction}
 \gamma_u=[\mathcal H_u,\id_{X_u},T_u],\qquad
 k_u\gamma_u=1_B,\qquad i_u\gamma_u=\delta_u.
\end{equation}

By \eqref{eq:sharp-exterior}, the class of $u$ is
\begin{equation}\label{eq:nz}
 n_z=d\boxtimes z+d^\sharp\boxtimes z^\sharp.
\end{equation}
Cycle invariance, additivity and Lemma~\ref{lem:contraction} now give
\begin{equation}\label{eq:recover-z}
 \mathcal B(n_z)=s_u\Delta_u,\qquad
 E_A\mathcal B(n_z)=\mathcal B(z).
\end{equation}
The retraction is at $k_u$, the corner of the shifted and doubled cycle.

\subsection{Fusion in the presented category}
\begin{theorem}\label{thm:dirac-fusion}
For $z\in\KK(A,B)$ and $y\in\KK(B,C)$ one has
\[
 \mathcal B(z)\mathcal B(y)=\mathcal B(zy)\qquad\text{in }\GKd.
\]
Hence $\mathcal A:\GKd\to\KK$ is an isomorphism of additive categories, with inverse $\mathcal B$. Every arrow has a simple form.
\end{theorem}
\begin{proof}
Use the cycle $u$ above. In analytic $KK$, put $v=\gamma_u y\in\KK(X_u,C)$ and choose a cycle representing $v$, hence a simple arrow $\mathcal B(v)$. Since $k_u$ is a homomorphism, \eqref{eq:precomp} and \eqref{eq:doubled-retraction} give
\[
 k_u\mathcal B(v)=\mathcal B(k_uv)=\mathcal B(y).
\]
All equalities in the following calculation take place in the presented category:
\begin{align*}
 \mathcal B(z)\mathcal B(y)
 &=E_As_u\Delta_u\mathcal B(y)
 &&\text{by \eqref{eq:recover-z}}\\
 &=E_As_u\Delta_uk_u\mathcal B(v)\\
 &=E_As_ui_u\mathcal B(v)
 &&\text{by the KK relation}\\
 &=E_A\mathcal B(s_ui_uv)
 &&\text{by \eqref{eq:precomp}}\\
 &=E_A\mathcal B(n_zy).
\end{align*}
The last equality uses the analytic identity
$s_ui_uv=s_ui_u\gamma_u y=s_u\delta_u y=n_zy$ and cycle invariance. Analytic exterior-product compatibility gives
\[
 n_zy=d\boxtimes(zy)+d^\sharp\boxtimes(z^\sharp y).
\]
Additivity and Lemma~\ref{lem:contraction} consequently imply
\[
 \mathcal B(z)\mathcal B(y)
 =E_A\mathcal B(d\boxtimes(zy))
  +E_A\mathcal B(d^\sharp\boxtimes(z^\sharp y))
 =\mathcal B(zy).
\]
The second branch contains $z^\sharp y$; the contraction lemma applies to this input directly.

This proves multiplicativity. Homomorphisms and synthetic generators are already simple forms, and no generators have been added. Hence every sum or composition equals $\mathcal B$ of its analytic image. Together with $\mathcal A\mathcal B=\id$, this proves $\mathcal B\mathcal A=\id$ and the asserted isomorphism.
\end{proof}

After the theorem, (T1) extends to an additive endofunctor on the quotient by
\[
 T(a)=\mathcal B(1_S\boxtimes\mathcal A(a)).
\]
Functoriality follows from the proved multiplicativity of $\mathcal B$. Every arrow is $\mathcal B(x)$ for its analytic image $x$, so (T2) gives $E_AT(a)=aE_B$ for every $a:A\to B$.

\section{An explicit obstruction inside the prescribed algebra X}\label{sec:obstruction}
In this and the next section we show the limits of the axiomatic fusion scheme of Sections~\ref{sec:sym} and~\ref{sec5} for graded algebras.

We now take $G$ trivial and $B=\C$. Let
\[
 A=C(S^1)\widehat\otimes\Cl,\qquad
 E=\ell^2(\mathbb Z)\otimes\C^{1|1}.
\]
The circle generator $u$ is even and the Clifford generator $\varepsilon$ odd. Write $Se_n=e_{n+1}$ and $Qe_n=e_n$ for $n\geq0$, $Qe_n=-e_n$ for $n<0$. Use
\[
 \sigma_x=\begin{pmatrix}0&1\\1&0\end{pmatrix},\quad
 \sigma_y=\begin{pmatrix}0&-i\\i&0\end{pmatrix},\quad
 \sigma_z=\begin{pmatrix}1&0\\0&-1\end{pmatrix}.
\]
The grading is $1\otimes\sigma_z$, and the cycle is
\begin{equation}\label{eq:example}
 \pi(u)=S\otimes1,\qquad c=\pi(\varepsilon)=1\otimes\sigma_y,
 \qquad F=Q\otimes\sigma_x.
\end{equation}
Indeed $F=F^*$, $F^2=1$, $Fc+cF=0$, and $[F,\pi(u)]=[Q,S]\otimes\sigma_x$ has finite rank. Density gives all Kasparov conditions. Adjoin only $(\C,0,0)$ and form $X_z$ as in \eqref{eq:XY}.

\subsection{The projection and partial isometry}
On $E$, direct multiplication gives
\[
 cF=\begin{pmatrix}-iQ&0\\0&iQ\end{pmatrix},\qquad
 q=\frac{\pi(1)+icF}{2}
 =\begin{pmatrix}P_+&0\\0&P_-\end{pmatrix},
\]
where $P_+=(1+Q)/2$ and $P_-=(1-Q)/2$. On the added $\C$ the operator is zero. Since $icF=-iF\pi(\varepsilon)$, we have $q\in X_z$. It is an even projection.

Set $v=q\pi(u)q$. Write $e_n^+,e_n^-$ for the even and odd basis vectors of $E$. On $qE$,
\[
 ve_n^+=e_{n+1}^+\ (n\geq0),\qquad
 ve_n^-=e_{n+1}^-\ (n\leq-2),\qquad ve_{-1}^-=0.
\]
Thus, for the even rank-one projections
\[
 e_+=\theta_{e_0^+,e_0^+},\qquad
 e_-=\theta_{e_{-1}^-,e_{-1}^-},
\]
we have
\begin{equation}\label{eq:defects}
 v^*v=q-e_-,\qquad vv^*=q-e_+.
\end{equation}
For an even projection $r\in X_z$, denote by $[r]\in KK(\C,X_z)$ the homomorphism class $\lambda\mapsto\lambda r$. Since $v$ is even, \eqref{eq:defects} and additivity on orthogonal projections imply
\begin{equation}\label{eq:equal}
 [e_-]=[e_+].
\end{equation}

\subsection{Where the opposite signs enter}
Both $e_+$ and $e_-$ are even projections, but their ranges are respectively an even and an odd line. The graded Morita equivalence $\K(E\oplus\C)\sim\C$ therefore sends their classes to $+1$ and $-1$.

The added corner is an even line. Compact even partial isometries identify any two even rank-one corners; graded Morita equivalence identifies an odd line with the negative of an even line in $KK$. Pushing these identities into $X_z$ gives
\begin{equation}\label{eq:opposite}
 [e_+]=[k_z],\qquad [e_-]=-[k_z].
\end{equation}
Together with \eqref{eq:equal}, this proves:
\begin{theorem}\label{thm:counter}
For the cycle \eqref{eq:example}, $2[k_z]=0$ in $KK(\C,X_z)$. There is no $\gamma\in KK(X_z,\C)$ with $k_z\gamma=m\,1_\C$ for any nonzero integer $m$.
\end{theorem}
\begin{proof}
Such a class would give $0=(2k_z)\gamma=2m$ in $KK(\C,\C)=\mathbb Z$.
\end{proof}

\section{The obstruction survives replacement of the cycle}\label{sec:invariant}
\begin{theorem}\label{thm:invariant}
For any cycle $z'$ representing the class of \eqref{eq:example}, padded with a single zero-representation $\C$-summand, one has
\[
 2[k_{z'}]=0\quad\text{in }KK(\C,X_{z'}).
\]
Thus no equivalent representative supplies the retraction required in Lemma~\ref{thm:lift}. In particular this class has no symmetry-admitting representative, even in the sense of Proposition~\ref{prop:mod}.
\end{theorem}
\begin{proof}
Write $z'=(E',\pi',F')$ and let $P_0=\pi'(1)$. Since $A$ is unital, compression of $F'$ to $E_0=P_0E'$ changes $F'\pi'(a)$ only by compacts. The complementary summand has zero representation and can be assigned operator zero. Put $c=\pi'(\varepsilon)|_{E_0}$. On $E_0$, $c$ is an odd selfadjoint unitary.

Replace the compressed operator by its selfadjoint part $F_s$ and then by $(F_s-cF_sc)/2$. These are compact changes: the defects are compact already after multiplying by $\pi'(1)$. Choose an even unitary identification $E_0=H\oplus H^{\op}$ under which
\[
 c=1\otimes\sigma_y,\qquad
 \pi'(u)=U\otimes1,\qquad
 F_a=Q\otimes\sigma_x.
\]
Here $U$ is unitary, $Q=Q^*$, $Q^2-1$ is compact and $[Q,U]$ is compact. These forms follow respectively from the odd unitarity of $c$, the exact commutation of $\pi'(u)$ with $c$, and the exact anticommutation of $F_a$ with $c$. Functional calculus changes $Q$ compactly into a selfadjoint involution; any finite-dimensional zero eigenspace can be assigned a sign. Denote the resulting operator on the full module by $\widetilde F$.

All these operations leave the concrete algebra unchanged:
\[
 (\widetilde F-F')\pi'(a)\in\K,\qquad
 \K+\pi'(A)+\widetilde F\pi'(A)=X_{z'}.
\]
No new Hilbert module is introduced in this normalization; we retain the full compact ideal $\K(E'\oplus\C)$, including the compact operators between all summands.

Set $P=(1+Q)/2$. The integer
\[
 n([z'])=\Ind(PUP:PH\to PH)
\]
is the odd Fredholm index pairing with the circle unitary after the standard Clifford identification \cite[Section~17.8]{Blackadar}
\[
 KK(C(S^1)\widehat\otimes\Cl,\C)\cong KK^1(C(S^1),\C).
\]
This pairing is a pairing of $KK$-classes, hence invariant under unitary equivalence, degenerate summands and Kasparov homotopy. We use the sign convention specified by the displayed compression formula. In the original example $PUP$ is the unilateral shift, so $n([z'])=n([z])=-1$.

Inside the unchanged $X_{z'}$ is the even projection
\[
 q=\frac{P_0+ic\widetilde F}{2}
 =\begin{pmatrix}P&0\\0&1-P\end{pmatrix}
\]
on $E_0$, extended by zero. The operator $a=q\pi'(u)q$ is even and Fredholm on $qE_0$. Its parity components have indices $n$ and $-n$, because their direct sum differs compactly from the unitary $U$. Its graded index is consequently $2n$.

We need the polar partial isometry \emph{inside this same algebra}, in order to obtain a relation in $KK(\C,X_{z'})$. Fredholmness gives closed range and a finite-dimensional kernel; thus $|a|$ has a spectral gap above zero. On $\Spec(|a|)$ the function
\[
 f_0(t)=\begin{cases}0,&t=0,\\ t^{-1},&t>0\end{cases}
\]
is continuous. Continuous functional calculus in the unital corner $qX_{z'}q$ therefore gives $v=a f_0(|a|)\in qX_{z'}q$. Its finite-rank even defect projections satisfy
\[
 v^*v=q-e_{\ker},\qquad vv^*=q-e_{\mathrm{coker}}.
\]
The even partial isometry $v$ equates these two projection classes in $KK(\C,X_{z'})$. Subtracting them yields
\[
 [e_{\ker}]-[e_{\mathrm{coker}}]=0
 \quad\text{in }KK(\C,X_{z'}).
\]
In the compact ideal the same difference corresponds, under graded Morita equivalence, to the graded index $2n$. The positive generator of that ideal maps to $[k_{z'}]$, because the added corner is one-dimensional and even. Hence
\begin{equation}\label{eq:invariant}
 2n([z'])[k_{z'}]=0.
\end{equation}
For $n=-1$ this proves the assertion. The final statement follows from Theorem~\ref{thm:sym} and Proposition~\ref{prop:mod}.
\end{proof}


\clearpage
\begin{thebibliography}{99}
\addcontentsline{toc}{section}{References}
\small
\bibitem{Blackadar} B. Blackadar, \emph{K-Theory for Operator Algebras}, second edition, MSRI Publications 5, Cambridge University Press, 1998.
\bibitem{Burgstaller} B. Burgstaller, Aspects of equivariant KK-theory in its generators and relations picture, \emph{Adv. Oper. Theory} \textbf{10} (2025), no.~2, article~32, 30 pp. \href{https://doi.org/10.1007/s43036-024-00412-y}{doi:10.1007/s43036-024-00412-y}.
\bibitem{flawedpaper} B. Burgstaller, \emph{The universal property of graded $KK^G$-theory}, earlier version, \href{https://arxiv.org/abs/2603.23157v1}{arXiv:2603.23157v1} (2026); subsequently withdrawn.
\bibitem{Cuntz} J. Cuntz, A new look at KK-theory, \emph{K-Theory} \textbf{1} (1987), 31--51.
\bibitem{Haag} U. Haag, Some algebraic features of $\mathbb Z_2$-graded KK-theory, \emph{K-Theory} \textbf{13} (1998), no.~1, 81--108.
\bibitem{HaagExt} U. Haag, On $\mathbb Z/2\mathbb Z$-graded KK-theory and its relation with the graded Ext-functor, \emph{J. Operator Theory} \textbf{42} (1999), no.~1, 3--36.
\bibitem{Higson} N. Higson, A characterization of KK-theory, \emph{Pacific J. Math.} \textbf{126} (1987), 253--276.
\bibitem{Kasparov81} G. G. Kasparov, The operator K-functor and extensions of $C^*$-algebras, \emph{Math. USSR-Izv.} \textbf{16} (1981), 513--572.
\bibitem{Kasparov88} G. G. Kasparov, Equivariant KK-theory and the Novikov conjecture, \emph{Invent. Math.} \textbf{91} (1988), 147--201.
\bibitem{Kucerovsky} D. Kucerovsky, The KK-product of unbounded modules, \emph{K-Theory} \textbf{11} (1997), 17--34.
\bibitem{LeGall} P.-Y. Le Gall, Th\'eorie de Kasparov \'equivariante et groupo\"ides. I, \emph{K-Theory} \textbf{16} (1999), 361--390.
\bibitem{MeyerEquivariant} R. Meyer, Equivariant Kasparov theory and generalized homomorphisms, \emph{K-Theory} \textbf{21} (2000), no.~3, 201--228. \href{https://doi.org/10.1023/A:1026536332122}{doi:10.1023/A:1026536332122}.
\bibitem{Meyer} R. Meyer, Categorical aspects of bivariant K-theory, in \emph{K-Theory and Noncommutative Geometry}, EMS Series of Congress Reports, European Mathematical Society, 2008, 1--39. \href{https://doi.org/10.4171/060-1/1}{doi:10.4171/060-1/1}.
\bibitem{vandenDungen} K. van den Dungen, Localisations of half-closed modules and the unbounded Kasparov product, \emph{Int. Math. Res. Not. IMRN} \textbf{2023}, no.~9, 7578--7615. \href{https://doi.org/10.1093/imrn/rnac065}{doi:10.1093/imrn/rnac065}.
\end{thebibliography}
\end{document}